\documentstyle[twoside,epsfig]{article}
\input epsf
\title{A note on some summations due to Ramanujan, their generalization and some allied series}
\author 
{A. K. Rathie\footnote{School of Mathematical and Physical Sciences, Central University of Kerala, Riverside Transit Campus, Padenenakad P. O. Nileshwar,
Kasaragad 671328, Kerala, India.
E-Mail: akrathie@rediffmail.com} 
\ and R. B. Paris\footnote{University of Abertay Dundee, Dundee DD1 1HG, UK.
E-Mail: r.paris@abertay.ac.uk}\\}
\begin{document}
\def\f#1#2{\mbox{${\textstyle \frac{#1}{#2}}$}}
\def\dfrac#1#2{\displaystyle{\frac{#1}{#2}}}
\def\boldal{\mbox{\boldmath $\alpha$}}
\newcommand{\bee}{\begin{equation}}
\newcommand{\ee}{\end{equation}}
\newcommand{\lam}{\lambda}
\newcommand{\ka}{\kappa}
\newcommand{\al}{\alpha}
\newcommand{\th}{\theta}
\newcommand{\om}{\omega}
\newcommand{\Om}{\Omega}
\newcommand{\fr}{\frac{1}{2}}
\newcommand{\fs}{\f{1}{2}}
\newcommand{\g}{\Gamma}
\newcommand{\br}{\biggr}
\newcommand{\bl}{\biggl}
\newcommand{\ra}{\rightarrow}
\newcommand{\mbint}{\frac{1}{2\pi i}\int_{c-\infty i}^{c+\infty i}}
\newcommand{\mbcint}{\frac{1}{2\pi i}\int_C}
\newcommand{\mboint}{\frac{1}{2\pi i}\int_{-\infty i}^{\infty i}}
\newcommand{\gtwid}{\raisebox{-.8ex}{\mbox{$\stackrel{\textstyle >}{\sim}$}}}
\newcommand{\ltwid}{\raisebox{-.8ex}{\mbox{$\stackrel{\textstyle <}{\sim}$}}}
\renewcommand{\topfraction}{0.9}
\renewcommand{\bottomfraction}{0.9}
\renewcommand{\textfraction}{0.05}
\newcommand{\mcol}{\multicolumn}
\date{}
\maketitle
\begin{abstract}
In this short note, we aim to discuss some summations due to Ramanujan, their generalizations and some allied series.

\vspace{0.4cm}

\vspace{0.3cm}

\end{abstract}

\vspace{0.3cm}

\begin{center}
{\bf 1. \  Introduction}
\end{center}
\setcounter{section}{1}
\setcounter{equation}{0}
\renewcommand{\theequation}{\arabic{section}.\arabic{equation}}
We start with the following summations due to Ramanujan \cite{R2}
\bee\label{e1}
1+\frac{1}{5}\left(\frac{1}{2}\right)^2+\frac{1}{9}\left(\frac{1\cdot3}{2\cdot 4}\right)^2+\cdots = \frac{\pi^2}{4\g^2(\f{3}{4})}
\ee
and
\bee\label{e2}
1+\frac{1}{5^2}\left(\frac{1}{2}\right)+\frac{1}{9^2}\left(\frac{1\cdot3}{2\cdot 4}\right)+\cdots = \frac{\pi^{5/2}}{8\surd 2\g^2(\f{3}{4})}.
\ee
As pointed out by Berndt \cite{B} the above summations can be obtained quite simply by putting (i) $a=b=\fs$, $c=\f{1}{4}$ and (ii) $a=\fs$, $b=c=\f{1}{4}$ in Dixon's summation theorem \cite[p. 52]{S} for the ${}_3F_2$ series, viz.
\[{}_3F_2\left[\!\!\begin{array}{c} a,\,b,\,c\\1+a-b,\,1+a-c\end{array}\!;1\right]=
\frac{\g(1+\fs a) \g(1+a-b) \g(1+a-c) \g(1+\fs a-b-c)}{\g(1+a) \g(1+\fs a-b) \g(1+\fs a-c) \g(1+a-b-c)}\]
valid provided Re\,$(\fs a-b-c)>-1$.

A similar series evaluation 
\bee\label{e3}
1+\frac{1}{5}\left(\frac{1}{2}\right)+\frac{1}{9}\left(\frac{1\cdot3}{2\cdot 4}\right)+\cdots = \frac{\pi^{3/2}}{2\surd 2\g^2(\f{3}{4})}
\ee
was also obtained by Ramanujan \cite{R1} using an integral representation. However, a more direct approach makes use of the fact that this series can be expressed as ${}_2F_1(\fs, \f{1}{4}; \f{5}{4}; 1)$ combined with the
well-known Gauss summation theorem \cite[p. 28]{S}
\bee\label{e3a}
{}_2F_1\left[\!\!\begin{array}{c}a,\,b\\c\end{array}\!;1\right]=\frac{\g(c) \g(c-a-b)}{\g(c-a) \g(c-b)}
\ee
valid provided Re\,$(c-a-b)>0$, by setting $a=\fs$, $b=\f{1}{4}$ and $c=\f{5}{4}$.

Next, let us consider the series
\bee\label{e4}
S=\frac{1}{b}+\frac{1}{2}\,\frac{1}{b+\mu}+\left(\frac{1\cdot 3}{2\cdot 4}\right)\,\frac{1}{b+2\mu}+\cdots ,
\ee
where\footnote{This can be extended to complex values of $\mu$ provided $|\arg\,\mu|<\pi$.} $\mu>0$. Then $S$ can be written as
\[S=\sum_{n=0}^\infty \frac{(\fs)_n}{n!}\,\frac{1}{b+n\mu}=\frac{1}{b}\sum_{n=0}^\infty \frac{(\fs)_n}{n!}\,\frac{(b/\mu)_n}{(b/\mu+1)_n}\]
\[=\frac{1}{b}\,{}_2F_1\left(\frac{1}{2},\, \frac{b}{\mu}; \,\frac{b}{\mu}+1; 1\right)\]
which can again be evaluated with the help of the Gauss summation theorem to yield the sum
\bee\label{e5}
S=\frac{\pi^{1/2}\,\g(b/\mu)}{\mu\,\g(b/\mu+\fs)}.
\ee
The case $\mu=2$ was considered by Ramanujan \cite{R1} using an integral representation. Clearly, the series (\ref{e4}) reduces to (\ref{e3}) by taking $\mu=4$, $b=1$. Thus the series (\ref{e4}) may be regarded as a generalization of (\ref{e3}).
\vspace{0.6cm}

\begin{center}
{\bf 2. \   Generalizations and other allied series}
\end{center}
\setcounter{section}{2}
\setcounter{equation}{0}
\renewcommand{\theequation}{\arabic{section}.\arabic{equation}}
In this section we shall mention some generalizations of Ramanujan's summations and also consider some allied series. For this, we apply the following results \cite{M, MP1}
\bee\label{e6}
{}_3F_2\left[\!\!\begin{array}{c} a,\, b, \,c\\b+m,\, c+1\end{array}\!;1\right]
=\frac{c \g(1-a) (b)_m}{(b-c)_m}\left\{\frac{\g(c)}{\g(1+c-a)}-\frac{\g(b)}{\g(1+b-a)} \sum_{k=0}^{m-1} \frac{(1-a)_k (b-c)_k}{(1+b-a)_k k!}\right\}
\ee
for positive integer $m$, and the generalized Karlsson-Minton summation formula for positive integers $(m_r)$ \cite{MP2, MS}
\bee\label{e9b}
{}_{r+2}F_{r+1}\left[\!\!\begin{array}{c}a, b,
\\c,\end{array}\!\!\!\!\begin{array}{c}(f_r+m_r)\\(f_r)
\end{array}\!;1\right]=
\frac{\g(c) \g(c-a-b)}{\g(c-a) \g(c-b)}\,\sum_{k=0}^m\frac{(-)^k(a)_k (b)_k C_k(r)}{(1+a+b-c)_k}
\ee
provided Re$\,(c-a-b)>m$, where $m=m_1+\cdots +m_r$, $(f_r)$ denotes the parameter sequence $(f_1, \ldots ,f_r)$ and the coefficients $C_k(r)$
are given by
\[C_k(r)=\frac{(-1)^k}{k!}\,{}_{r+1}F_r\left[\!\!\begin{array}{c}-k,\\{}\end{array}
\!\!\!\begin{array}{c}(f_r+m_r)\\(f_r)\end{array}\!;1\right].\]
In the particular case $r=1$, Vandermonde's summation theorem can be used to show that
\[C_k(1)=\left(\!\!\!\begin{array}{c}m\\k\end{array}\!\!\!\right)\,\frac{1}{(f)_k}.\]

Let us consider the extension of (\ref{e4}) in the form
\bee\label{e7}
S=\sum_{n=0}^\infty \frac{(\fs)_n}{n!}\,\frac{1}{(1+n/b)(1+n/c)}=
{}_3F_2\left[\!\!\begin{array}{c}\fs,\,b,\,c\\b+1,\,c+1\end{array}\!;1\right].
\ee
This can be evaluated by means of (\ref{e6}) with $m=1$ to give
\[S=\left\{\begin{array}{ll}\dfrac{\pi^{1/2} bc}{b-c} \left\{\dfrac{\g(c)}{\g(c+\fs)}-\dfrac{\g(b)}{\g(b+\fs)}\right\} & (b\neq c)\\
\\
\dfrac{\pi^{1/2}b^2 \g(b)}{\g(b+\fs)} \{\psi(b+\fs)-\psi(b)\} & (b=c).\end{array}\right.
\]
The special case $b=c$ is obtained by a limiting process with $\psi$ denoting the logarithmic derivative of the gamma function. If we let $b=c=\f{1}{4}$ in (\ref{e7}), we immediately obtain Ramanujan's summation (\ref{e2}), upon noting that $\psi(\f{3}{4})-\psi(\f{1}{4})=\pi$.

Further, consider the series
\bee\label{e8}
S_p=\sum_{n=0}^\infty \left(\frac{(\fs)_n}{n!}\right)^{\!2}\frac{1}{(n+1) \ldots (n+p)},
\ee
where integer $p\geq 1$. This corresponds to
\[S_p=\frac{1}{p!}\,{}_2F_1(\fs, \fs; p+1; 1)\]
which upon use of Gauss' summation theorem (\ref{e3a}) reduces to
\bee\label{e9a}
S_p=\frac{\g(p)}{\g^2(p+\fs)}.
\ee
Thus, from (\ref{e8}) and (\ref{e9a}) it follows that
\[S_1=\frac{4}{\pi},\qquad S_2=\frac{16}{9\pi}, \qquad S_3=\frac{128}{225\pi}\]
and so on.

If we let $a=b=\fs$ and $c=p+1$ for positive integer $p$ in (\ref{e9b}), then we have when $r=1$ (with $m=m_1$)
\[{}_3F_2\left[\!\!\begin{array}{c}\fs,\,\fs,\ f+m\\p+1,\,f\end{array}\!;1\right]=\frac{p! \g(p)}{\g^2(p+\fs)} \sum_{k=0}^m (-)^k\left(\!\!\!\begin{array}{c}m\\k\end{array}\!\!\!\right)\frac{((\fs)_k)^2}{(f)_k (1-p)_k}\qquad(p>m).\]
When $m=1$, we therefore find
\bee\label{e9}
\sum_{n=0}^\infty \left(\frac{(\fs)_n}{n!}\right)^{\!2}\frac{n+f}{(n+1) \ldots (n+p)}=
\frac{\g(p)}{\g^2(p+\fs)}\left(f+\frac{1}{4(p-1)}\right)
\ee
for $p=2, 3, \ldots$, and when $m=2$
\bee\label{e10}
\sum_{n=0}^\infty \left(\frac{(\fs)_n}{n!}\right)^{\!2}\frac{(n+f)(n+f+1)}{(n+1) \ldots (n+p)}=
\frac{\g(p)}{\g^2(p+\fs)}\left(f(f+1)+\frac{f+1}{2(p-1)}+\frac{9}{16(p-1)(p-2)}\right)
\ee
for $p=3, 4, \ldots\,$.

When $r=2$ and $m_1=m_2=1$ (so that $m=2$), we obtain from (\ref{e9b})
\[{}_4F_3\left[\!\!\begin{array}{c}\fs,\,\fs,\ f_1+1,\,f_2+1\\p+1,\,f_1,\,f_2\end{array}\!;1\right]=\frac{p!}{f_1 f_2}\sum_{n=0}^\infty\left(\frac{(\fs)_n}{n!}\right)^{\!2}\frac{(n+f_1)(n+f_2)}{(n+1) \ldots (n+p)}\]
\[=\frac{p! \g(p)}{\g^2(p+\fs)}\left\{1+\frac{C_1(2)}{4(p-1)}+\frac{9C_2(2)}{16(p-1)(p-2)}\right\},\]
where $C_1(2)=(1+f_1+f_2)/(f_1f_2)$ and $C_2(2)=1/(f_1 f_2)$. Hence
\bee\label{e11}
 \sum_{n=0}^\infty\left(\frac{(\fs)_n}{n!}\right)^{\!2}\frac{(n+f_1)(n+f_2)}{(n+1) \ldots (n+p)}=\frac{\g(p)}{\g^2(p+\fs)}\left(f_1 f_2+\frac{f_1+f_2+1}{4(p-1)}+\frac{9}{16(p-1)(p-2)}\right)
\ee
for $p=3, 4, \ldots\,$.

We remark that series such as (\ref{e9}) can also be obtained by a `telescoping'
process applied to the series $S_p$ in (\ref{e8}). For example, it easy to see that
 \[\sum_{n=0}^\infty \left(\frac{(\fs)_n}{n!}\right)^{\!2}\frac{n+f}{(n+1)\ldots (n+p)}=
S_{p-1}+(f-p)S_p.\]

\end{document}